\documentclass[preprint,12pt]{elsarticle}
\usepackage{amsmath}
\usepackage{amsfonts}
\usepackage{bm}
\usepackage{lipsum}
\usepackage{enumerate}
\usepackage{graphicx}
\usepackage{amsfonts}

\usepackage{mathrsfs}
\usepackage{subfigure}
\usepackage{epstopdf}
\usepackage{url}
\usepackage{verbatim}
\usepackage{amssymb}
\usepackage{hyperref}
\usepackage{color}
\usepackage{listings}

\def\@begintheorem#1#2{\par\bgroup{\scshape #1\ #2. }\it\ignorespaces}
\def\@opargbegintheorem#1#2#3{\par\bgroup%
   {\scshape #1\ #2\ ({\upshape #3}). }\it\ignorespaces}
\def\@endtheorem{\egroup}

\if@onethmnum
  \newtheorem{theorem}{Theorem}
  \newtheorem{lemma}[theorem]{Lemma}
  \newtheorem{corollary}[theorem]{Corollary}
  \newtheorem{proposition}[theorem]{Proposition}
  \newtheorem{definition}[theorem]{Definition}
\newtheorem{example}[theorem]{Example}
\newtheorem{remark}[theorem]{Remark}
\newtheorem{homework}[theorem]{Homework}
\newtheorem{case}[theorem]{}

\else

\fi

\usepackage[T1]{fontenc}

\journal{XXX }

\begin{document}

\begin{frontmatter}
\title{Tauberian identities and the connection to Wile E. Coyote physics}

\author[1]{Roberto Camassa}
\ead{camassa@amath.unc.edu}
\author[1]{Richard M. McLaughlin \corref{mycorrespondingauthor}}
\ead{rmm@email.unc.edu}
\cortext[mycorrespondingauthor]{Corresponding author} 
\address[1]{Department of Mathematics, University of North Carolina, Chapel Hill, NC, 27599, United States}

\begin{abstract}
 The application of the motion of a vertically suspended mass-spring system released under tension is studied focusing upon the delay timescale for the bottom mass as a function of the spring constants and masses.  This ``hang-time", reminiscent of the Coyote and Road Runner cartoons, is quantified using the far-field asymptotic expansion of the bottom mass' Laplace transform.  These asymptotics are connected to the short time mass dynamics through Tauberian identities and explicit residue calculations.  It is shown, perhaps paradoxically, that this delay timescale is maximized in the large mass limit of the top ``boulder".  Experiments are presented and compared with the theoretical predictions.  This system is an exciting example for the teaching of mass-spring dynamics in classes on Ordinary Differential Equations, and does not require any normal mode calculations for these predictions.  
\end{abstract}

\begin{keyword}
Road Runner \sep Wile E. Coyote \sep Hangtime \sep Tauberian Identities 
\end{keyword}
\end{frontmatter}

\section{Introduction}

The evolution of mass-spring systems is a subject which has received great attention for its modeling of many different physical and biological systems, for its connection in limiting cases to wave dynamics, and for simple examples demonstrating the normal mode decomposition in many texbooks on Ordinary Differential Equations \cite{Boyce}.

Here, we examine a less studied and exciting example of the dynamics of a hanging mass-spring system initially in dynamic equilibrium which is suddenly released from its top suspended point and allowed to evolve due to the transient spring forces and the gravitational force, assumed to be Hookean in the former and constant in the latter.  

This interesting dynamics we will see presents behavior reminiscent of Wile E. Coyote of the famous Road-Runner cartoons, in which Wile E. is frequently tricked by the Road Runner into falling off a cliff, often with a large boulder falling above him.  Wile E. Coyote frequently, in these circumstances, appears to pause momentarily, with no clear force explaining the delay in his falling, before plummeting to the canyon floor.  Our dynamics, experimentally, computationally, and analytically, demonstrate a similar delay timescale.  We quantify this timescale analytically as a function of the number and amount of masses as well as the values of the spring constants.  

This study presents possibly the simplest physical example demonstrating the power of Tauberian identities relating the short time asymptotics of a function with the far field expansion of its Laplace Transform.  As such it makes an excellent example for the study of coupled systems in any level class on ordinary differential equations.   

In section 2, we describe the experiment and discuss the experimental observations.  In section 3, we set up the mathematical model for the linear coupled system of $n$ masses connected by $n-1$ springs, initially hung from above, and suddenly released from its hanging point.  We compute the initial mass distribution in equilibrium, and then compute the short time asymptotics of the bottom mass which we show is very flat, scaling as $C_n t^{2n}$.  In section 4 we present the quantitative comparison of the experimental dynamics with the mass-spring model dynamics, focussing upon the asymptotic prediction of the delay timescale in predicting the experimental hang-time.  In section 5, we discuss these results as well as mathematical differences associated in the continuum limit and compare and contrast with common similarities and misconceptions regarding the finite and infinite propagation of information in partial differential equations, such as in the wave and heat equations.

\section{Wile E. Coyote Physics:  falling stretched mass-spring systems }
\begin{figure}
\centering
\includegraphics[width=0.75\linewidth]{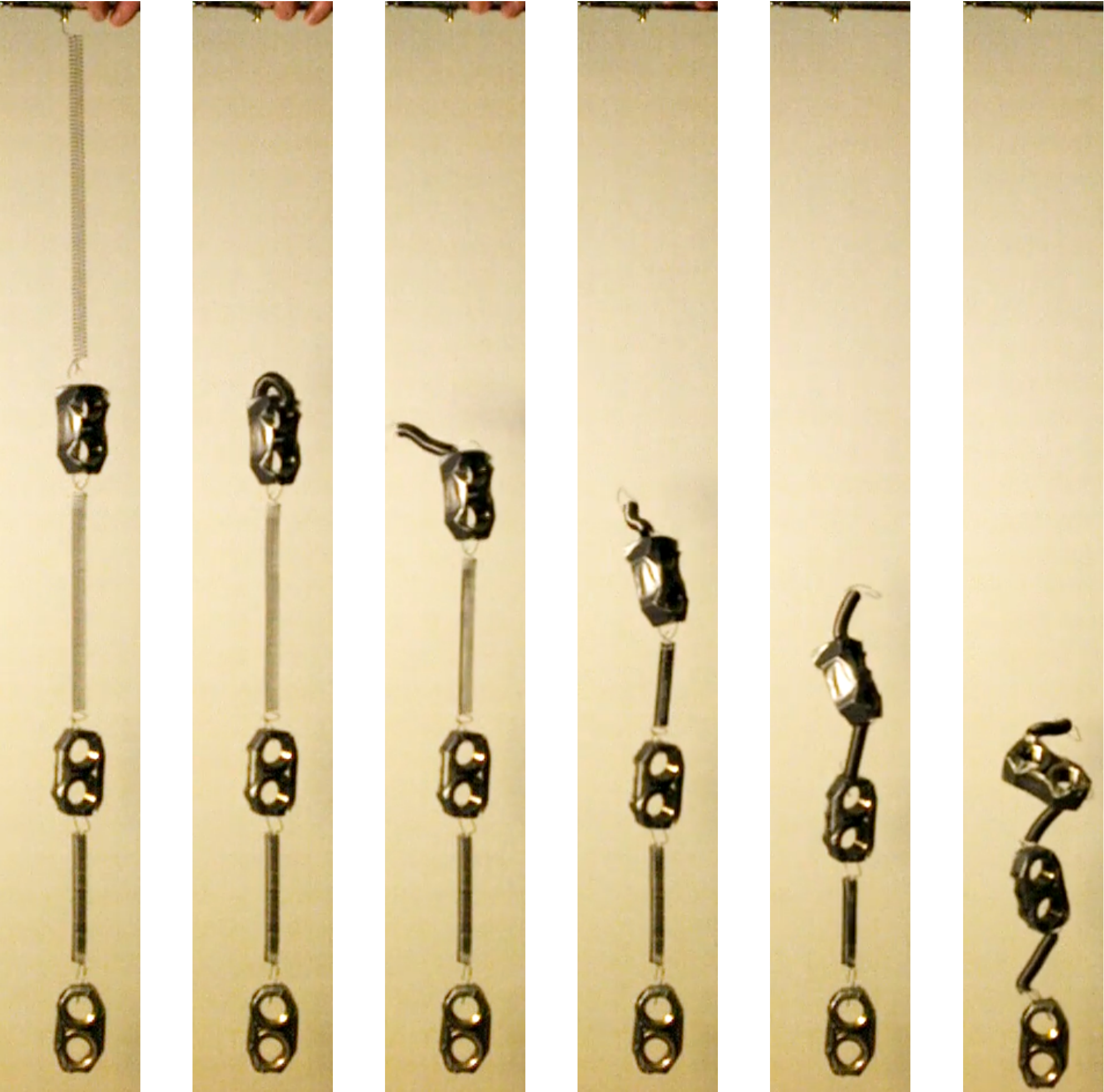}
\includegraphics[width=0.3\linewidth]{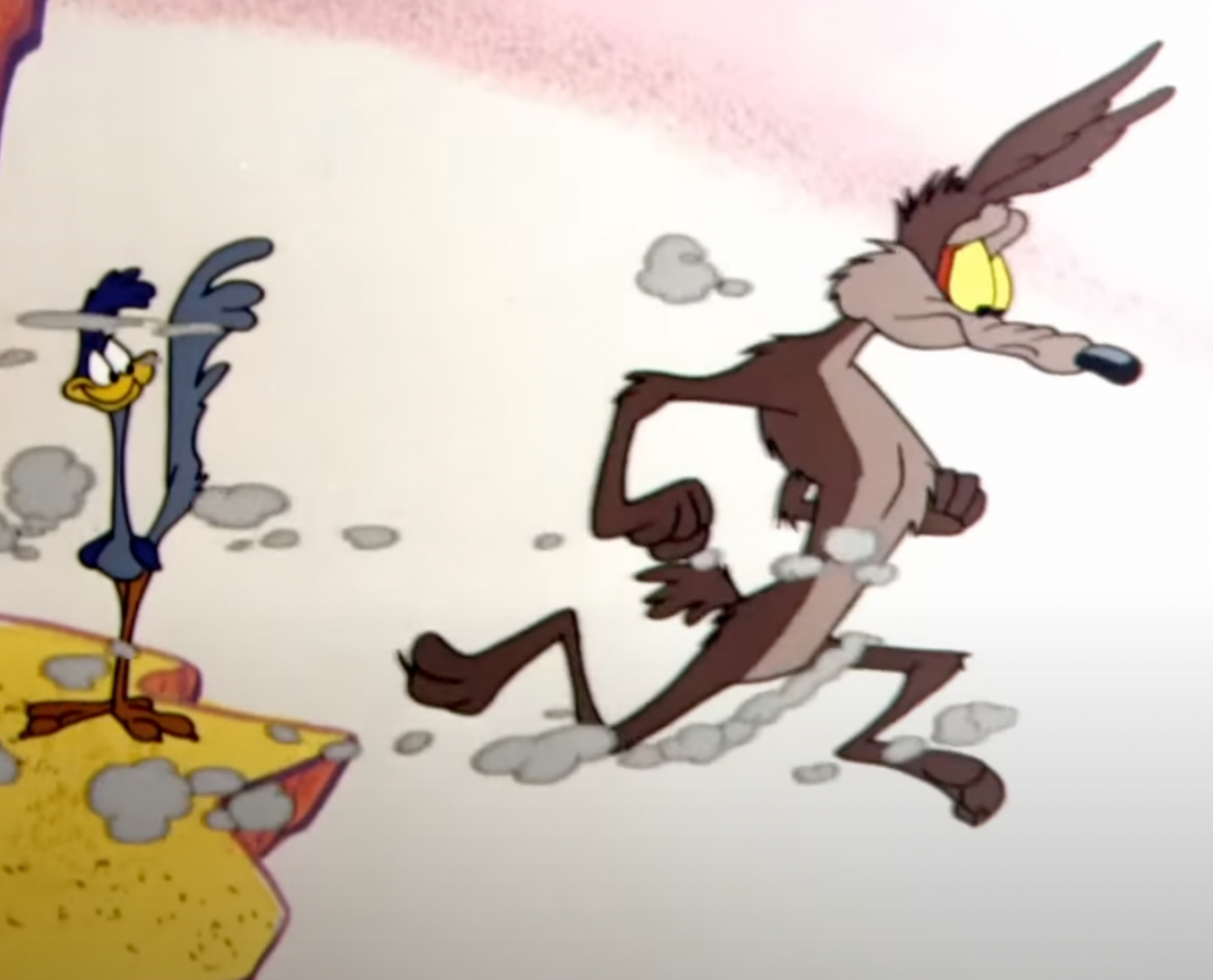}
\includegraphics[width=0.19\linewidth]{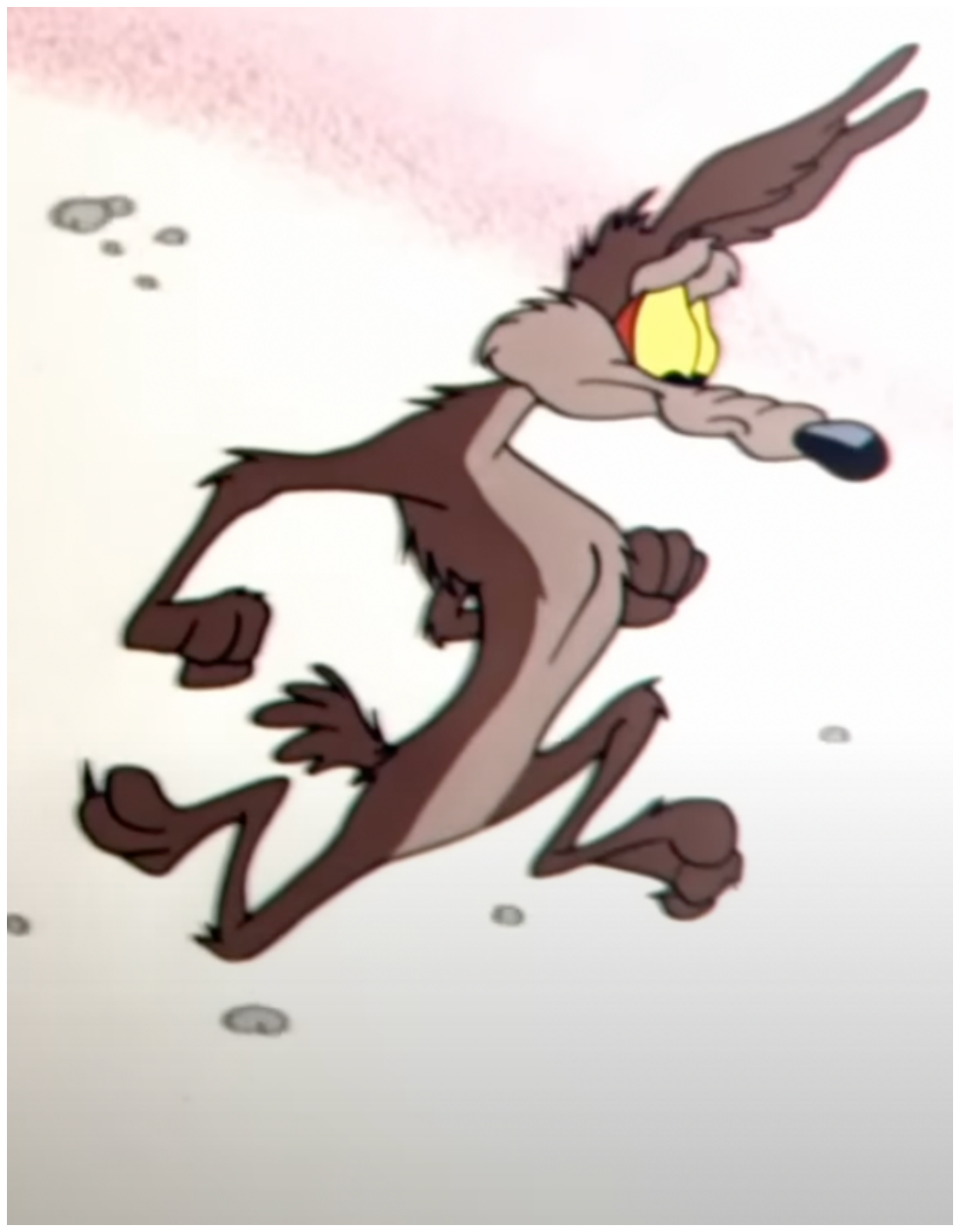}
\includegraphics[width=0.14\linewidth]{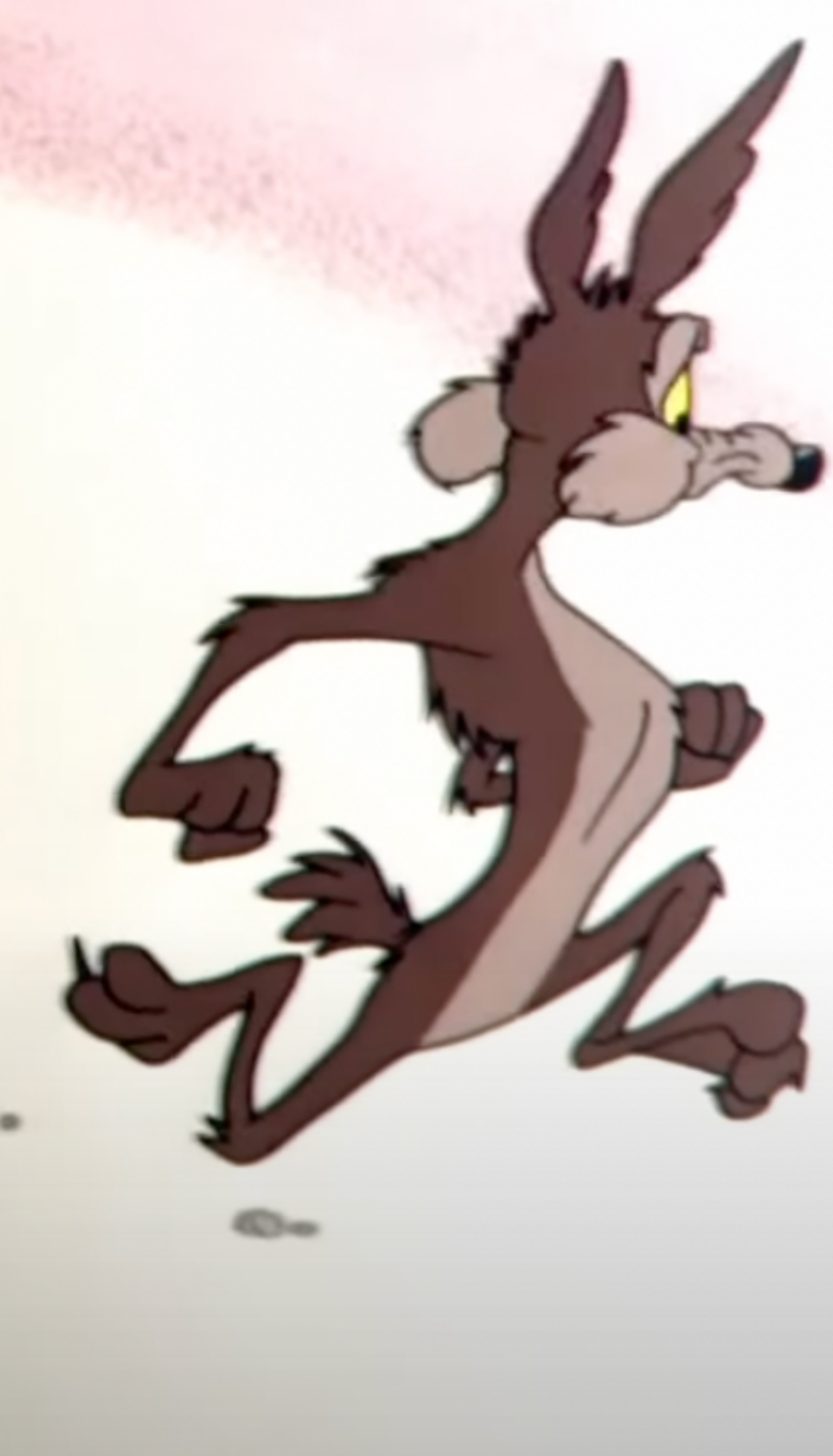}
\includegraphics[width=0.14\linewidth]{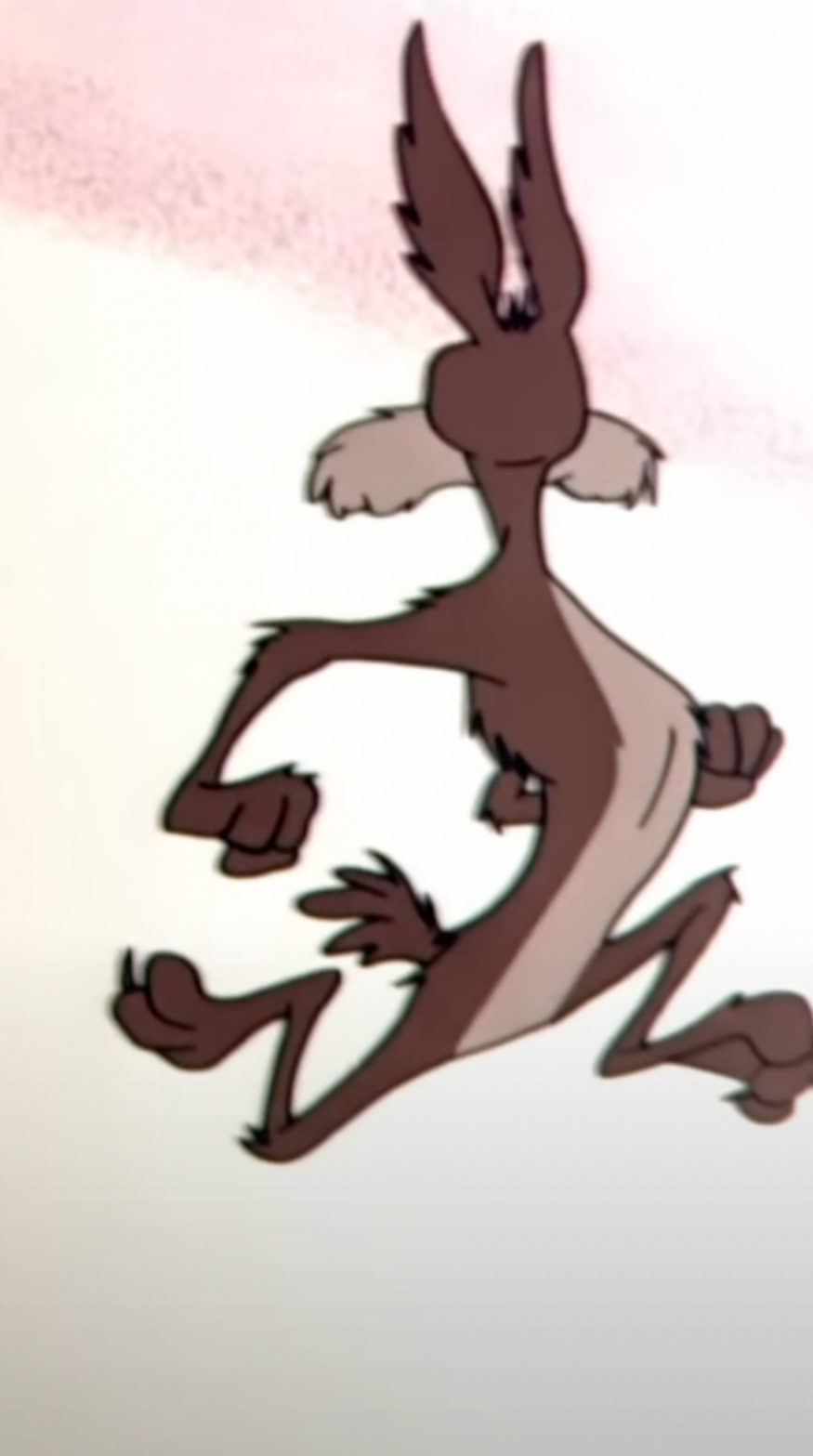}
\includegraphics[width=0.11\linewidth]{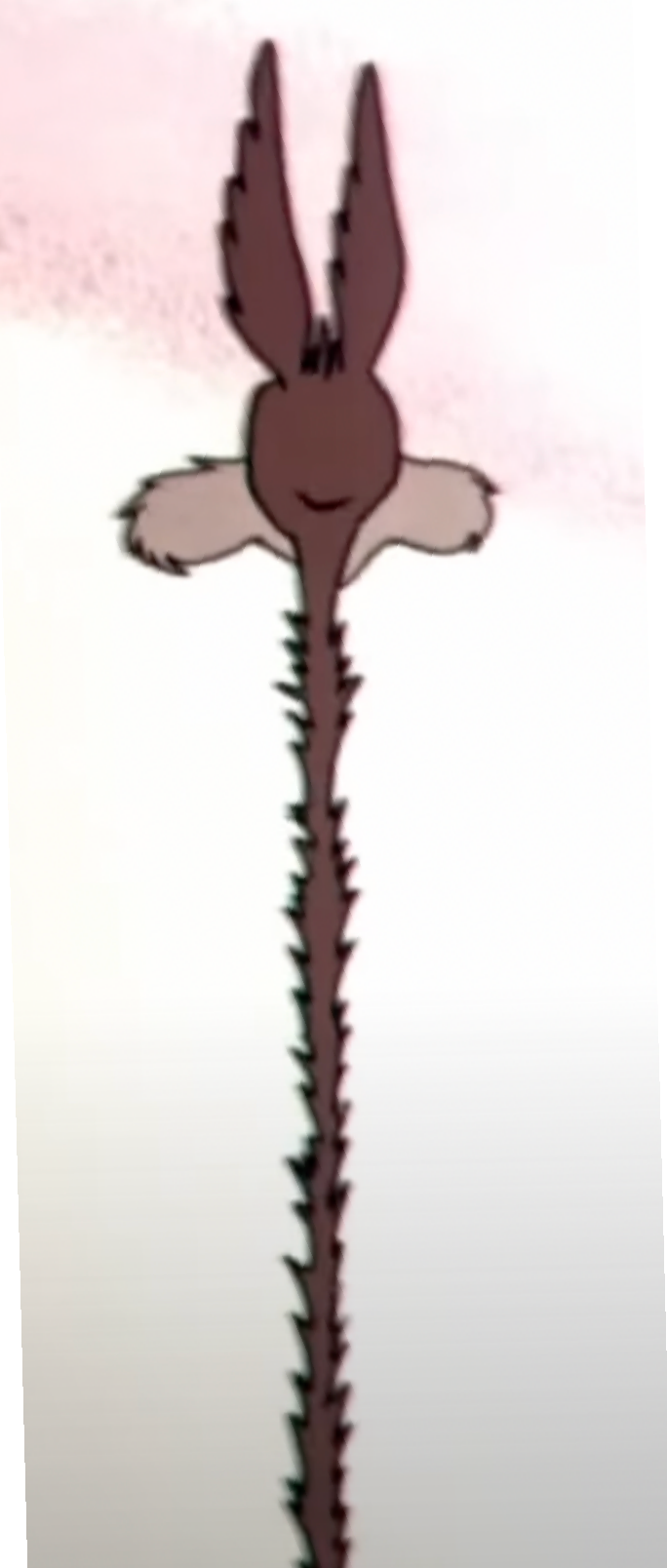}
\caption{Falling mass-spring system:  Top:  initially in static equilibrium, the top spring is released, and the subsequent evolution is observed on $0.03$ second intervals.  No explicit top mass is affixed to the top spring here.  Bottom:  sample Wile E. Coyote montage depicting a similar delay timescale before falling.}
\label{montage}
\end{figure}

We next outline our experiments with 4 masses attached via 3 linear Hooke-law springs.  Images of the system are captured with and Edgertronic high speed digital video camera at a frame rate of 1000 frames per second.  Figure \ref{montage} shows time series of the behavior taken at equal time slices after the top mass is released:  in the top panel in this trial, there is no mass attached to the top spring, the first panel shows the hanging system in static equilibrium just before releasing the top spring. Subsequent panels are given at uniform time intervals of $0.03$s. Last, in the bottom panel is a sample Wile E. Coyote delay from one of the famous cartoons. See accompanying videos for the full evolution.  Several things are worth observing. First, there is essentially no motion of any of the bottom three masses until an observable time delay passes.  After approximately $0.17$ seconds, the bottom mass begins to fall.  

\begin{figure}
\centering
\includegraphics[width=1.\linewidth]{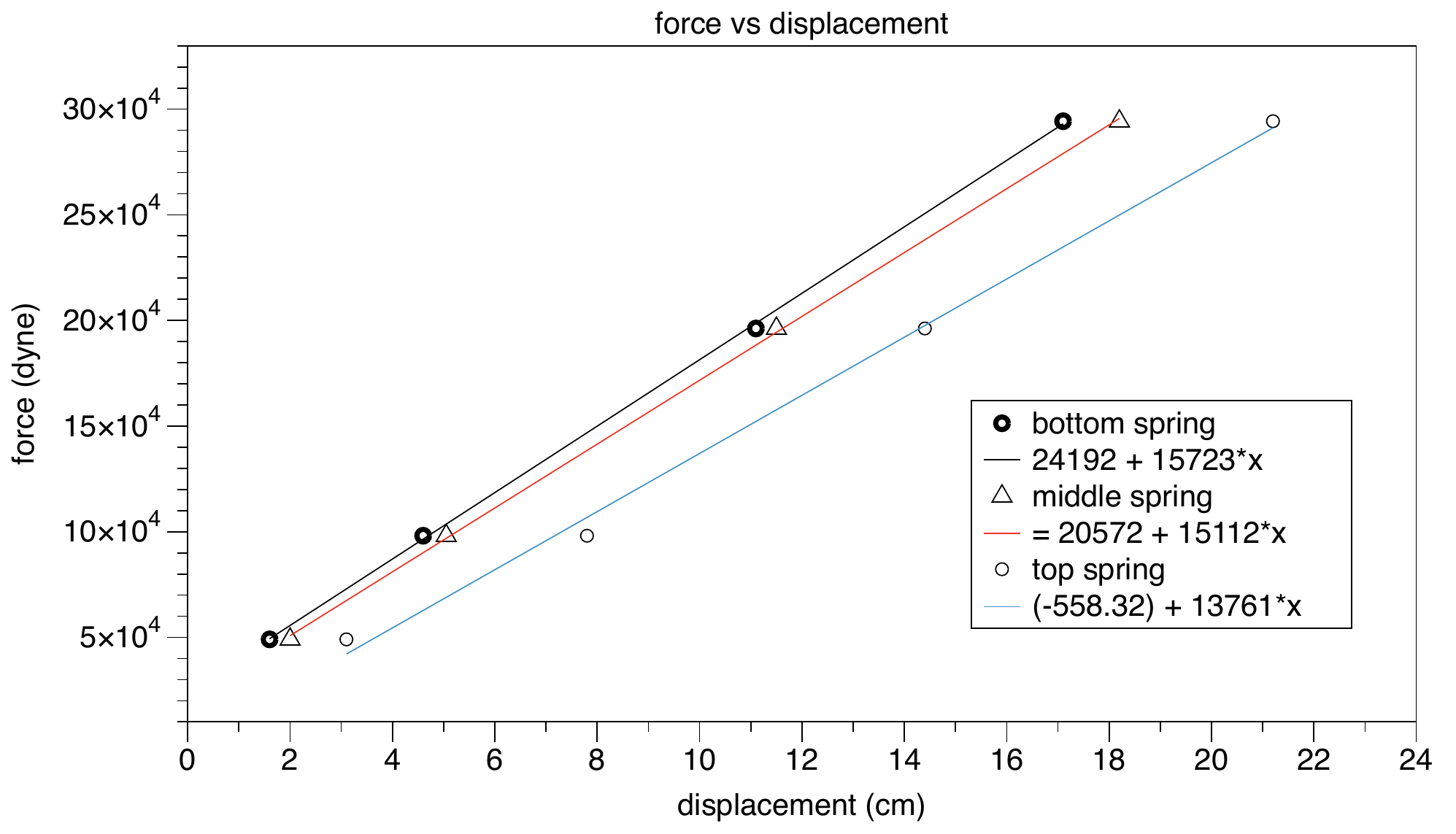}
\caption{Spring constants:  The three springs equilibrium positions are measured at static equilibrium under an array of different mass loadings to measure the spring constants.  Linear fits are provided and the slopes provide the spring constants, measured in dynes/cm.}
\label{springconstant}
\end{figure}

The springs were wound in house with a Porter Spring Winders tool.  In figure \ref{springconstant}, we present the measurement of the Hookean spring constants for each of the three springs.  This is done by imaging the static equilibrium displacement of each spring under various mass loadings, and fitting the force vs displacement data to a straight line.  The slope of which provides the spring constants, measured in dynes/cm.  The top spring has $k_1=13,761$ dynes/cm, middle with $k_2=15,112$ dynes/cm, and bottom with $k_3=15,723$ dynes/cm.  We note that while the springs were all wound equally, the top and middle springs were subjected to greater loading during the experiment, and consequently were stretched and weakened comparatively with the bottom spring.  The mass of each spring is $4.44$ grams.  The bottom mass is $m_4=107.28$ grams, third mass is $m_3=106.56$ grams, and the second mass is $m_2=107.5$ grams.  In this trial, we do not attach a mass to the top of the top spring,  $m_1=0$ grams.  Since the system is being released, the motion of the system will depend upon the spring masses (which do not play an immediate role in the normal mode calculations for studying the dynamics of the hanging system).  As such, the spring masses are each added to masses $m_1$, $m_2$, and $m_3$ in the analysis below.  

\section{III. Mathematical modeling with linear coupled mass-spring system}
Assuming that our $n-1$ spring constants are linear, with Hookean spring constants $k_j$, as is clear from the data presented in Figure \ref{springconstant}, the dynamics of the suspended system can be written as 
\begin{eqnarray}
m_1 \frac{d^2 x_1}{dt^2}&=& k_1(x_2-x_1)-m_1 g +N \nonumber\\ 
m_j \frac{d^2 x_j}{dt^2}&=& k_{j-1}(x_{j-1}- x_j) +k_j(x_{j+1}-x_j) -m_j g, \qquad j=2,\cdots,n-1\\
m_n \frac{d^2 x_n}{dt^2}&=& k_{j-1}(x_{n-1}-x_n) -m_{n} g\nonumber
\end{eqnarray}
Here, the normal force, $N=\sum_{j=1}^n m_j g$, is the force required to hold the mass system in place so that it is initially in equilibrium.  While we see that it is not necessary, the equilibrium distribution of positions, $(x_1(0),x_2(0),\cdots,x_n(0))$ may be determined by solving the linear system
\begin{eqnarray}
\begin{pmatrix}
k_1 & -k_2 & 0& 0 &\cdots& 0\\
-k_1 & k_1+k_2 & -k_2&0&\cdots&0\\
0&\ddots&\ddots&\ddots&0&0\\
\vdots&0&\ddots&\ddots&\ddots&0\\
0&\cdots&0&-k{n-2}&k_{n-2}+k_{n-1}&-k_{n-1}\\
0&0&\cdots&0&-k_{n-1}&k_{n-1}
\end{pmatrix}
\quad
\begin{pmatrix}
x_1(0)\\
x_2(0)\\
x_3(0)\\
\vdots\\
x_{n-1}(0)\\
x_n(0)
\end{pmatrix}
&=&
\begin{pmatrix}
(\sum_{j=2}^nm_j)g\\
m_2 g\\
m_3 g\\
\vdots\\
m_{n-1} g\\
m_n g
\end{pmatrix}
\end{eqnarray}
We note that this matrix is not invertible on account of the fact that the system has a one degree of freedom translation invariance arising from the observation that the spring lengths are irrelevant in the calculation.  This is overcome by setting $x_1(0)=0$, and solving the resulting system as an upper triangular system yielding a unique solution. This solution does not matter as it does not in any way affect the dynamics as we next see from changing from global to local coordinates for each mass:
\begin{eqnarray}
z_1(t)&=&x_1(t)\\
z_j(t)&=&x_j(t)-x_j(0),\qquad j=2,\cdots,n
\end{eqnarray}
The new system of equations then is 
\begin{eqnarray}
m_1 \frac{d^2 z_1}{dt^2}&=& k_1(z_2-z_1)+N-\sum_{j=1}^n m_j g \nonumber\\ 
m_j \frac{d^2 z_j}{dt^2}&=& k_{j-1}(z_{j-1}-z_j) +k_{j}(z_{j+1}-z_j) , \qquad j=2,\cdots,n-1 \\
m_n \frac{d^2 z_n}{dt^2}&=& k_{n-1}(z_{n-1}-z_n) \nonumber
\label{eqsystem}
\end{eqnarray}
Clearly, this system is in equilibrium as the normal force, $N$, cancels the sum in the first equation, assuming that all initial positions and velocities are zero, i.e., $z_j(0)=\frac{dz_j}{dt}|_{t=0}=0$ for all $j$. 

We next cut the normal force by simply setting $N=0$ in system \ref{eqsystem}, and then watching the subsequent system evolution.  This is equivalent to releasing the top mass, as done in the prior experimental section.  Before computing the ensuing dynamics, it is useful to non-dimensionalize this system.  Let the system time scale be set by the frequency of the first equation, $T=\sqrt{\frac{m1}{k_1}}$ and the system length scale, $L=T^2( \frac{\sum_{j=1}^n m_j g}{m_1})$.  Let $y_j=\frac{z_j}{L}, \tau=\frac{t}{T}$, then the non-dimensionalized system reads:
\begin{eqnarray}
\frac{d^2 y_1}{d\tau^2}&=& (y_2-y_1)-1 \nonumber\\  
\frac{d^2 y_j}{d\tau^2}&=& \alpha_{j} (y_{j-1}- y_j) +\beta_j(y_{j+1}-y_j) , \qquad j=2,\cdots,n \\
\frac{d^2 y_n}{d\tau^2}&=& \alpha_n(y_{n-1}-y_n) \nonumber
\end{eqnarray}
where the non-dimensional spring constants are $\alpha_j=\frac{m_1}{m_j}\frac{k_{j-1}}{k_1}$ and  $\beta_j=\frac{m_1}{m_j}\frac{k_{j}}{k_1}$ for $j=2,\cdots,n-1$.

To compute the ensuing evolution, we compute the Laplace transform of the system, $\hat y_j(s)=\int_0^\infty dt e^{-s \tau} y_j(\tau)$.  This results in the following matrix system:
\begin{eqnarray}
\begin{pmatrix}
s^2 +1 & -1 & 0& 0 &\cdots& 0\\
-\alpha_2 & s^2 + \alpha_2+\beta_2 & -\beta_2&0&\cdots&0\\
0&\ddots&\ddots&\ddots&0&0\\
\vdots&0&\ddots&\ddots&\ddots&0\\
0&\cdots&0&-\alpha_{n-1}&s^2 + \alpha_{n-1}+\beta_{n-1}&-\beta_{n-1}\\
0&0&\cdots&0&-\alpha_n&s^2 + \alpha_n
\end{pmatrix}
\quad
\begin{pmatrix}
\hat y_1(s)\\
\hat y_2(s)\\
\hat y_3(s)\\
\vdots\\
\hat y_{n-1}(s)\\
\hat y_n(s)
\end{pmatrix}
=
\begin{pmatrix}
\frac{-1}{s}\\
0\\
0\\
\vdots\\
0\\
0
\end{pmatrix}
\end{eqnarray}
Let matrix, $A$ be defined as
\begin{eqnarray}
A=
\begin{pmatrix}
s^2 +1 & -1 & 0& 0 &\cdots& 0\\
-\alpha_2 & s^2 + \alpha_2+\beta_2 & -\beta_2&0&\cdots&0\\
0&\ddots&\ddots&\ddots&0&0\\
\vdots&0&\ddots&\ddots&\ddots&0\\
0&\cdots&0&-\alpha_{n-1}&s^2 + \alpha_{n-1}+\beta_{n-1}&-\beta_{n-1}\\
0&0&\cdots&0&-\alpha_n&s^2 + \alpha_n
\end{pmatrix}
\end{eqnarray}
The fastest way to invert this system is using Cramer's rule.  The Laplace transform of the last mass' position, $\hat y_n(s)$, is then the ratio of two determinants:
\begin{eqnarray*}
\hat y_n(s)= \frac{\det B}{\det A}
\end{eqnarray*}
where the matrix B is 
\begin{eqnarray}
B&=& \begin{pmatrix}
s^2 +1 & -1 & 0& 0 &\cdots& \frac{-1}{s}\\
-\alpha_2 & s^2 + \alpha_2+\beta_2 & -\beta_2&0&\cdots&0\\
0&\ddots&\ddots&\ddots&0&0\\
\vdots&0&\ddots&\ddots&\ddots&0\\
0&\cdots&0&-\alpha_{n-1}&s^2 +\alpha_{n-1}+\beta_{n-1}&0\\
0&0&\cdots&0&-\alpha_n&0
\end{pmatrix}
\end{eqnarray}
Do to the sparse driver vector, the determinant of B is immediately evaluated as 
\begin{eqnarray*}
\det B&=& -\left(\prod_{j=2}^n\alpha_j\right)\frac{1}{s}
\end{eqnarray*}
We remark that similar formulas are available, and computable for any of the mass positions, which we present below.  

{\bf Short time asymptotics:} We also comment that here is where Tauberian asymptotics help to greatly detangle the determinant of $A$.   In general, the full inverse involves computing the complete determinant of $A$, and its inverse transform will involve a lengthy sum of residues which will yield the standard normal form decomposition of the system.  But it is important to consider the experiment:  no oscillations are observed in the dynamics, it is all about computing the short time delay before the last mass is detected to move.  To this end, computing the short time asymptotics will suffice in predicting the evolution.  To this end, the Tauberian identities dictate that the short time dynamics of the n'th mass position is set by the far field, large $s$ asymptotics of the Laplace transform.  To whit, we need only compute the large $s$ asymptotic expansion of the determinant of $A$ to compute the short time asymptotic expansion of the solution.  Clearly, this is set by the strong diagonal dominance of matrix $A$:
\begin{eqnarray}
\det A &=& P_{2n}(s)\sim \prod_{j=1}^n A_{j,j} \\
&=&  \prod_{j=1}^n (s^2 + \alpha_j+\beta_j) \\
&\sim& s^{2 n} \qquad \mbox{as} \quad s\rightarrow \infty
\end{eqnarray}
So we have that the short time asymptotics of the position, $y_n(\tau)$, of the last mass is given by the following Bromwich complex integral:
\begin{eqnarray}
y_n(\tau)&\sim& -\frac{(\prod_{j=2}^n \alpha_j)}{2\pi i} \int_{\gamma-i\infty}^{\gamma+i \infty}\frac{e^{s \tau}}{s^{2n+1}}ds\\
&=&-\frac{(\prod_{j=2}^n \alpha_j)}{2\pi i} \oint \frac{1+(s\tau) + \cdots \frac{(s \tau)^{2n}}{(2n)!}+\cdots}{s^{2n+1}}ds
\end{eqnarray}
Where the closed integral is taken around any counter-clockwise contour encircling the origin.  The result is: 
\begin{eqnarray}
y_n(\tau)&\sim&-\left(\frac{m_1^{n-1}\prod_{j=2}^{n-1} k_j}{k_1^{n-2} \prod_{j=2}^n m_j}\right) 
\frac{\tau^{2n}}{{(2n)}!}, \quad \mbox{as} \quad \tau\rightarrow 0^+
\label{taub}
\end{eqnarray}
In dimensional form, the result it:
\begin{eqnarray}
z_n(t)&\sim& -Q \frac{ t^{2n}}{(2n)!}, \quad \mbox{as} \quad t\rightarrow 0^+
\label{shorttime}
\end{eqnarray}
where 
\begin{eqnarray*}
Q&=& \left(\frac{(\prod_{j=1}^{n-1} k_j)\sum_{j=1}^n m_j g}{\prod_{j=1}^n m_j}\right) \end{eqnarray*}
This has a few interesting limits to discuss.  First, if the first mass (boulder) is very large, say compared to the other masses, $m_j=m$ for $j=2,\cdots,n$, and taking all the spring constants equal, $k_j=k$ for $j=1,\cdots, n-1$, then the prefactor reads $Q=(1+\frac{(n-1)m}{m_1}) \frac{k^{n-1}g}{m^{n-1}}$, which is monotonically decreasing with increasing $m_1$, with a non-zero limit for large $m_1$.  Consequently, holding all other parameters fixed, increasing $m_1$ will lead to an increased hang-time before the $n$'th particle starts moving, with the maximum delay occurring for infinite $m_1$ (see discussion below regarding computing the hang-time).  

We remark that the alternative limit of vanishing mass, $m_1$ is quite interesting to analyze in dimensional coordinates (or with an alternative non-dimensionalization), and provides a re-derivation of the short time power law dynamics, albeit with a different power law in this limit.  To this end, consider the system above in \ref{eqsystem}:  in the limit of small $m_1$, the first equation reads, $m_1\frac{d^2 z_1}{dt^2}= -k (z_1-z_2) - m_1g - \sum_{j=2}^n m_jg$.  Here we have extracted the first mass from the sum.  Clearly this equation is well approximated, for small $m_1$ by $z_1=z_2-\frac{\sum_{j=2}^n m_j g}{k}$.  Inserting this into the second equation, gives 
\begin{eqnarray*}
m_2\frac{d^2 z_2}{dt^2}&=& \frac{k}{m_2} (z_1-2z_2+z_3)\\
&=& k(z_3-z_2) -\sum_{j=2}^n m_j g
\end{eqnarray*}
Notice, by comparison with the system in (\ref{eqsystem}), the system is identical, with the first mass removed altogether.  This is equivalent to reducing the mass counter to $n-1$, and the result is the same in \ref{shorttime}, with $n$ replaced by $n-1$, and the $j$ index running from $2$ to $n$.  In this case, the short time behavior is less flat, scaling as $\tau^{2n-2}$.

We also point out that the so-called Tauberian identities in this example may be understood immediately by scaling.  Keeping the general determinant of the matrix, $A$ as a polynomial of degree $2n$, $P_{2n}(s)=s^{2n}+a_{2n-2} s^{2n-2}+\cdots+a_0$, the exact solution is represented as the contour integral:
\begin{eqnarray*}
y_n(\tau)&=&\frac{-(\prod_{j=2}^n \alpha_j)}{2\pi i} \oint_{\cal C} \frac{e^{\tau s}}{s P_{2n}(s)}ds
\end{eqnarray*}
Where the contour, $\cal{C}$ is a closed circle containing all the roots of the denominator.  We may rescale the integral through $z=\tau s$, the result is
\begin{eqnarray*}
\oint_{\cal C} \frac{e^{\tau s}}{s P_{2n}(s)}ds&=& \oint_{\cal \tilde C}\frac{e^z}{z P_{2n}(\frac{z}{\tau})}dz\\
&=&\tau^{2n}\oint_{\cal \tilde C}\frac{e^z}{z (z^{2n}+a_{2n-2}\tau^2 z^{2n-2} + \cdots  + a_0\tau^{2n} ) }dz\\
&\sim&  \tau^{2n}\oint_{\cal \tilde C}\frac{e^z}{z^{2n+1} }dz\qquad \tau\rightarrow 0^+
\end{eqnarray*}
Here, $\cal \tilde C$ is the rescaled contour by $\tau$.  Of course, by independence of path, this contour can be stretched arbitrarily.  This justifies the Tauberian identity used to arrive at \ref{taub}.

We next provide the result for the short time asymptotics of the $j$'th mass position.  To this end, the Tauberian asymptotics yields the far field asymptotics of the generalized matrix, $\tilde B$, with driver replacing the $j'$'th column of the matrix, $A$.  The result is
\begin{eqnarray*}
\det \tilde B &\sim& -\prod_{j=2}^{j}\alpha_j \qquad \tau \rightarrow 0^+
\end{eqnarray*}
Then the short time asymptotics of the $j$'th mass' position, $y_j(t)$ is, for $j=2,\cdots, n$: 
\begin{eqnarray*}
y_j(t)&\sim& - Q_j \frac{t^{2j}}{(2j)!}\quad \mbox{as} \quad t\rightarrow 0^+\\
Q_j&=&\left(\frac{(\prod_{i=1}^{j-1} k_i)\sum_{j=1}^n m_j g}{\prod_{i=1}^j m_i}\right) \end{eqnarray*}
We note that the first mass enjoys $y_1\sim -\frac{\sum_{j=1}^n (m_j g)}{2 m_1} t^2$ as $t \rightarrow 0^+$.  Interestingly, the first mass initially experiences an acceleration substantially larger than gravity on account of the energy stored in the stretched springs.  Of course intuition suggests that this acceleration cannot persist indefinitely, as long time asymptotics shows.  

{\bf Long time asymptotics:} First, note that the center of mass of the system does fall accelerating with the gravitational constant, $g$, which is immediately observed by summing the equations of motion.  One may ask if indeed the long time asymptotics of any individual mass also falls at with this same constant acceleration or not?  To that end, we immediately can compute the leading order, long time asymptotic expansion of the last mass, $y_n(\tau)$ using similar Tauberian identities, except in this case taking the limit of long time.  To that end, repeating the calculation instead at long time:
\begin{eqnarray*}
\oint_{\cal C} \frac{e^{\tau s}}{s P_{2n}(s)}ds&=& \oint_{\cal \tilde C}\frac{e^z}{z P_{2n}(\frac{z}{\tau})}dz\\
&=&\tau^{2n}\oint_{\cal \tilde C}\frac{e^z}{z (z^{2n}+a_{2n-2}\tau^2 z^{2n-2} + \cdots  + a_0\tau^{2n} ) }dz\\
&\sim& \oint_{\cal \tilde C}\frac{e^z}{z(a_0+a_2 z^2/\tau^2)}  dz\qquad \tau\rightarrow \infty
\end{eqnarray*}
It is an interesting question to ask how the coefficients of the characteristic polynomial of a matrix depend upon the matrix coefficients themselves.  This sometimes is described in terms of the so-called Newton identities.  Some results are known and immediate:  the lowest order coefficient, $a_0$ is proportional to the determinant of the matrix, and the second to highest power has $a_{2n-2}$ proportional to the trace of the matrix.  In this case, the matrix in question is $A$ with the parameter $s$ set to zero.  This matrix has determinant zero.  So the dominant asymptotics at long time reads:
\begin{eqnarray*}
y_n(\tau) &\sim&  \frac{-(\prod_{j=2}^n \alpha_j)}{2\pi i}\tau^2\oint_{\cal \tilde C}\frac{e^z}{z^3 a_2}dz\\
&=& -(\prod_{j=2}^n \alpha_j)\frac{\tau^2}{2 a_2} \qquad \tau\rightarrow \infty
\end{eqnarray*}
So the long time asymptotics are set by the coefficient, $a_2$.  Calculating this term requires some co-factor expansions, and the result is 
\begin{eqnarray*}
a_2=\frac{(\prod_{j=2}^n \alpha_j)}{m_1} \sum_{j=1}^n m_j 
\end{eqnarray*} 
leading to the result that $z_n(t)\sim -g t^2/2$ as $t \rightarrow \infty$.  This establishes that the $n$th mass falls with the acceleration of gravity at long time, and that the corrections are sub-dominant, meaning that there will not be large oscillations observable at long time.  Of course, the actual experiment does extend long enough to observe this behavior. 

Lastly, we comment on corrections to the short time asymptotics.  These are immediately provided by including the first and second highest powers of $s$ in formula for, $P_{2n}(s)$:

\begin{eqnarray*}
y_n(\tau)&=&\frac{-(\prod_{j=2}^n \alpha_j)}{2\pi i} \oint_{\cal C} \frac{e^{\tau s}}{s P_{2n}(s)}ds\\
&\sim& \frac{-(\prod_{j=2}^n \alpha_j)}{2\pi i}\tau^{2n}\oint_{\cal \tilde C}\frac{e^z}{z (z^{2n}+a_{2n-2}\tau^2 z^{2n-2} ) }dz,\qquad \tau\rightarrow 0^+\\
&\sim& \frac{-(\prod_{j=2}^n \alpha_j)}{2\pi i}\tau^{2n}\oint_{\cal \tilde C}\frac{e^z}{ z^{2n+1} } (1-a_{2n-2} \tau^2/z^2)dz\\
&=& -\frac{(\prod_{j=2}^n \alpha_j)\tau^{2n}}{(2n)!} \left(1- \frac{a_{2n-2} \tau^{2}}{(2n+2)(2n+1)}\right)
\end{eqnarray*}
Here, the coefficient $a_{2n-2}$ is equal to the trace of the matrix $A$ with $s=0$.  We note that in the analysis section below, this correction is very small for the timescales observed in the experiment, with the leading order asymptotics accurately predicting the hang-time for all the masses observed in our experimental mass-spring system.

\section{IV. Model Comparison with Experimental Observations}
Here we first discuss how this formula indeed sets a proper timescale for the hang-time in this system, and then proceed to demonstrating the predictive power of our asymptotic formulae at short time.  

While it is true that algebraic expression $\tau^{2n}$ is non-zero for any positive $\tau$, detectable information requires an instrument be able to measure the motion.  In our case, with our high speed Edgertronic camera, we can set an instrumentation threshold by the length scale of one pixel.  As such, we can define the hang-time, $t_h$ to be the amount of elapsed time between when the top mass is released and when the bottom mass moves a detectable distance in the camera image.  Having a high enough frame rate is necessary to resolve the times in between top mass release and the hang-time.  In our theory, the hang-time is given by 
\begin{eqnarray}
t_h&=& \left(\frac{L_p((2n)!)}{Q}\right)^{\frac{1}{2n}} \\
&=& \left(L_p(2n)!
\left(\frac{\prod_{j=1}^n m_j}{(\prod_{j=1}^{n-1} k_j)\sum_{j=1}^n m_j g}\right)\right)^{\frac{1}{2n}}
\end{eqnarray}
where $L_p$ is the detectable distance. 

\begin{figure}
\centering
\includegraphics[width=0.44\linewidth]{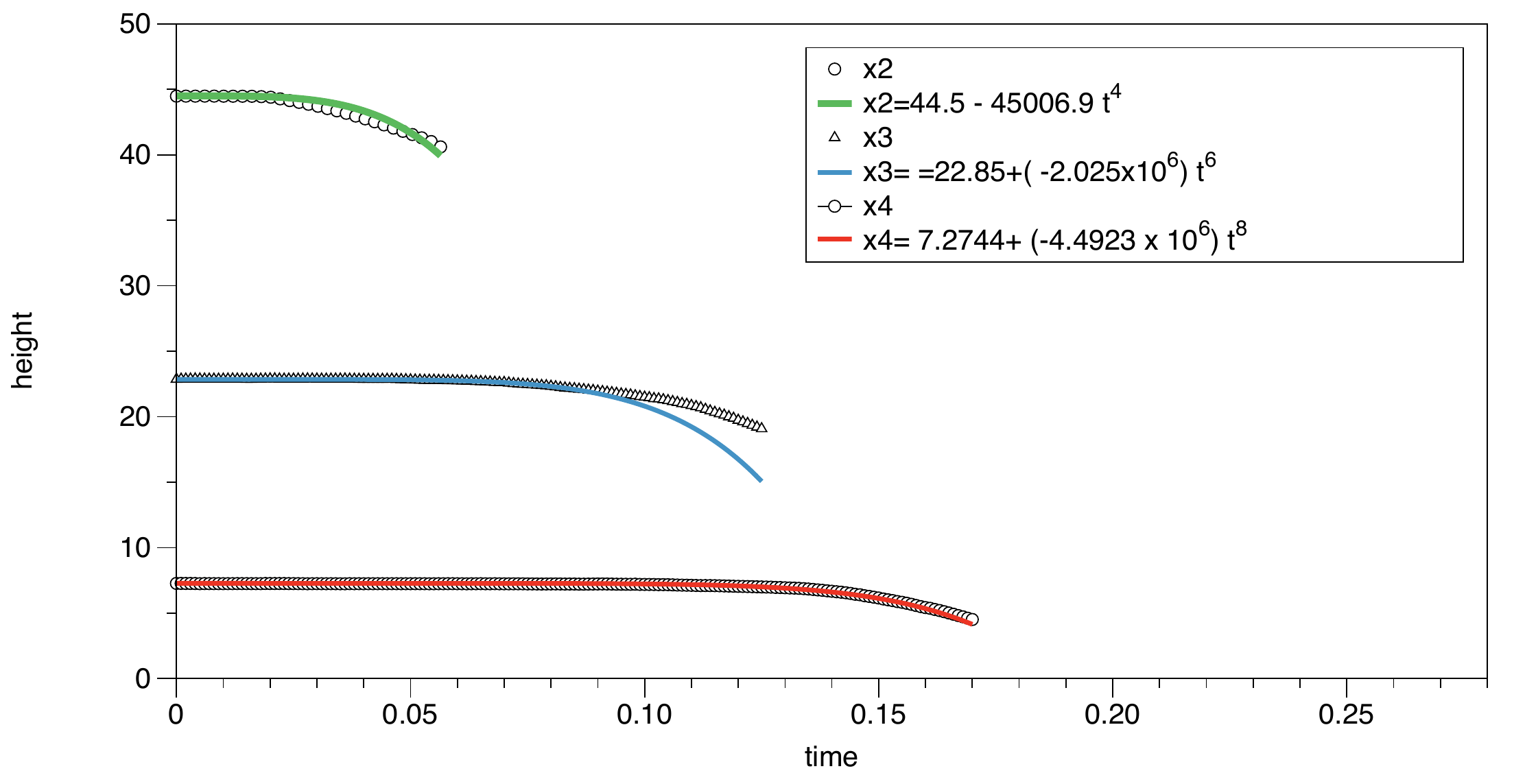}
\includegraphics[width=0.4\linewidth]{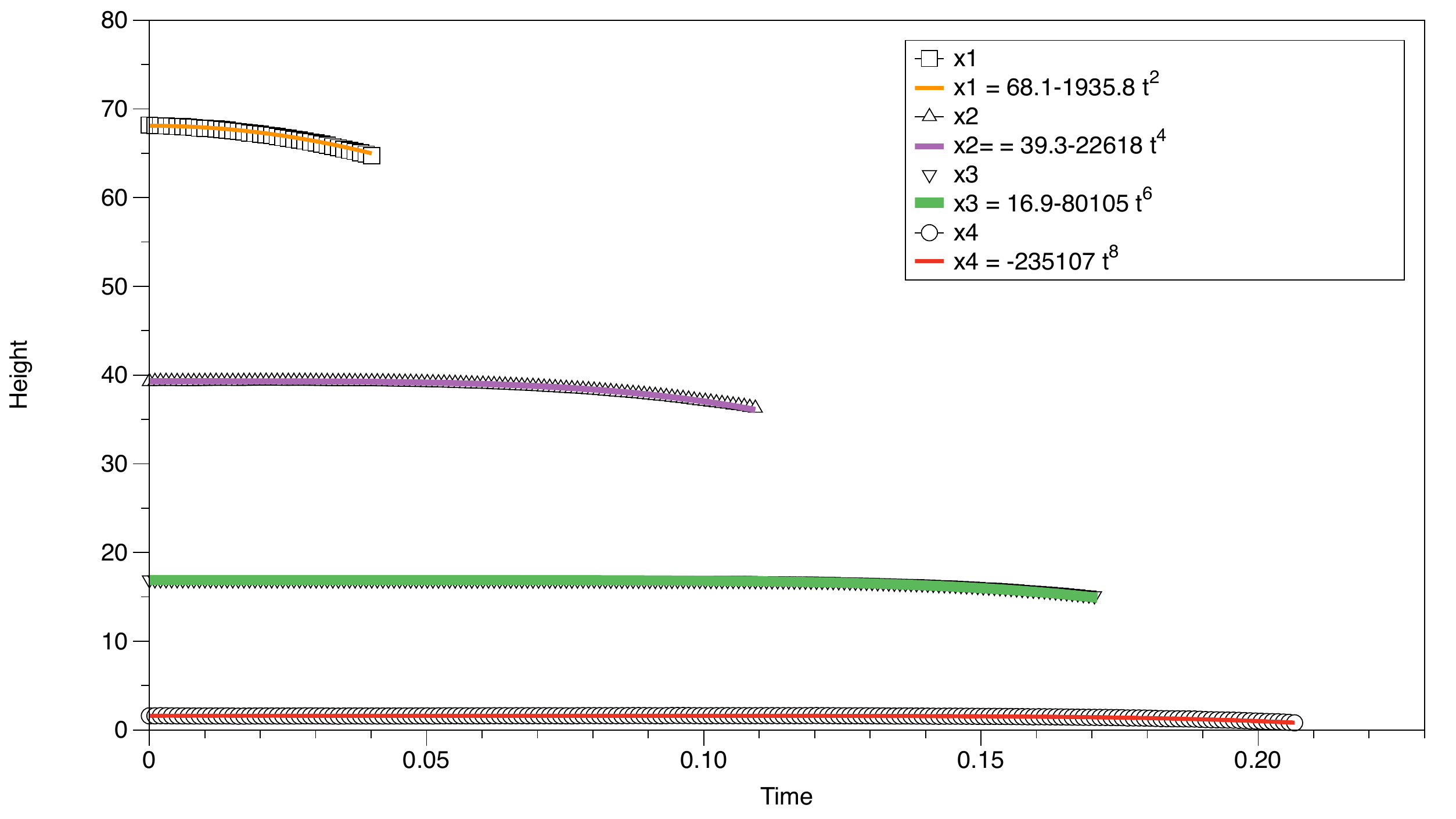}
\caption{Experimental mass tracks with leading order small time asymptotic theory superimposed with no adjustable parameters, and theoretical values reported in the legend in cgs units. Left case:  no explicit top mass is attached as in Figure 1.  Right panel, explicit top mass, $m_1=107.7$g.   }
\label{3tracks}
\end{figure}

Alternatively, we can simply plot the trajectories and superimpose the power law short time asytmptotics.  In figure \ref{3tracks}, we plot the short time asymptotic theory super-imposed upon the experimentally tracked mass heights for two different cases, the latter with no extra mass beyond the spring affixed to the top spring, and in the former with $m_1=107.7$ grams.  A few things to note here.  First, for the case on the left which has no explicit top mass, the theory does a fairly good job for all three masses, where we set $m_1$ to be equal to the top spring mass.  Clearly, the spring masses play a role in the subsequent dynamics through the inertia they provide, but the modeling is not expected to be as good as for the case where there is an explicit mass attached to the top spring.  That is the case in the right panel, and in that case all four masses motion are quite accurately predicted by these formulas.  Further note that the hang time of mass $m_4$ is considerably longer than that in the left panel on account of the larger top mass in the right panel.  This prediction is born out by our theory.  Second, incorporating the short asymptotic corrections does not provide much improvement on these timescales.  In fact, by truncating the geometric series to two terms as done above, gives rise to strong competition with the leading order term being negative, and the correction being positive, which gives rise to the erroneous prediction that the particle rises on longer time scales.  But this error is not observable on the short timescales of the experiments.  It is worth noting that summing the full Laurent expansion gives a correction which is negative definite for all time, but again over the timescales of the experiment, this calculation is indistinguishable from two term expansion presented above.  

A few comments on the experiments.  First, the lighting in the first trial with no explicit top mass was superior to the second experiment.  Second, the second trial exhibited some minor out of plane motion of the bottom mass, and additionally the poor lighting made tracking more difficulty.  Nonetheless, the modeling is far better in this case, again on account of the fact that there is a well-defined top mass in this trial.  The trackings were performed in the Image Tank program, using one of the color (RGB) fields, and tracking a contour value.

\begin{figure}
\centering
\includegraphics[width=0.47\linewidth]{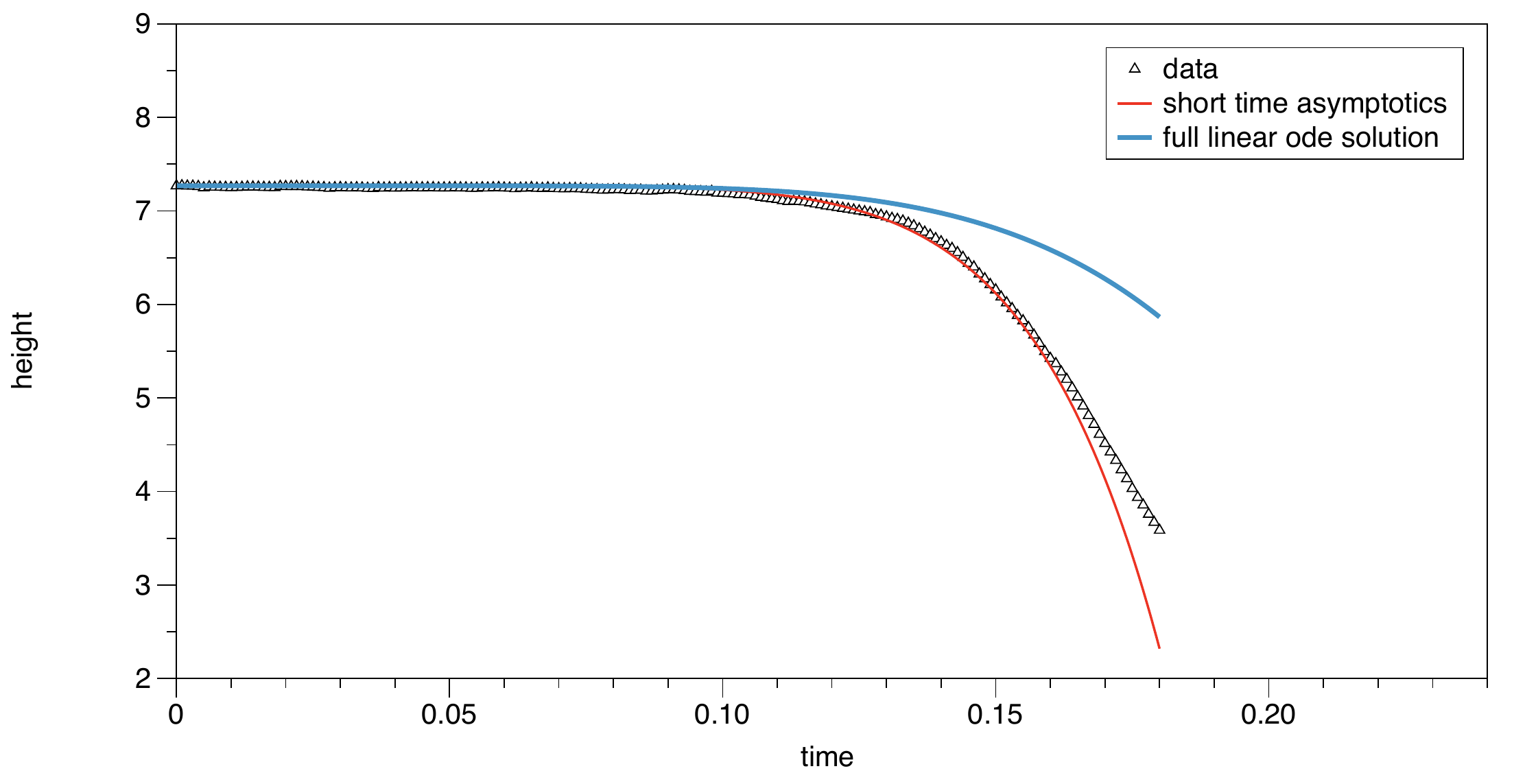}
\includegraphics[width=0.43\linewidth]{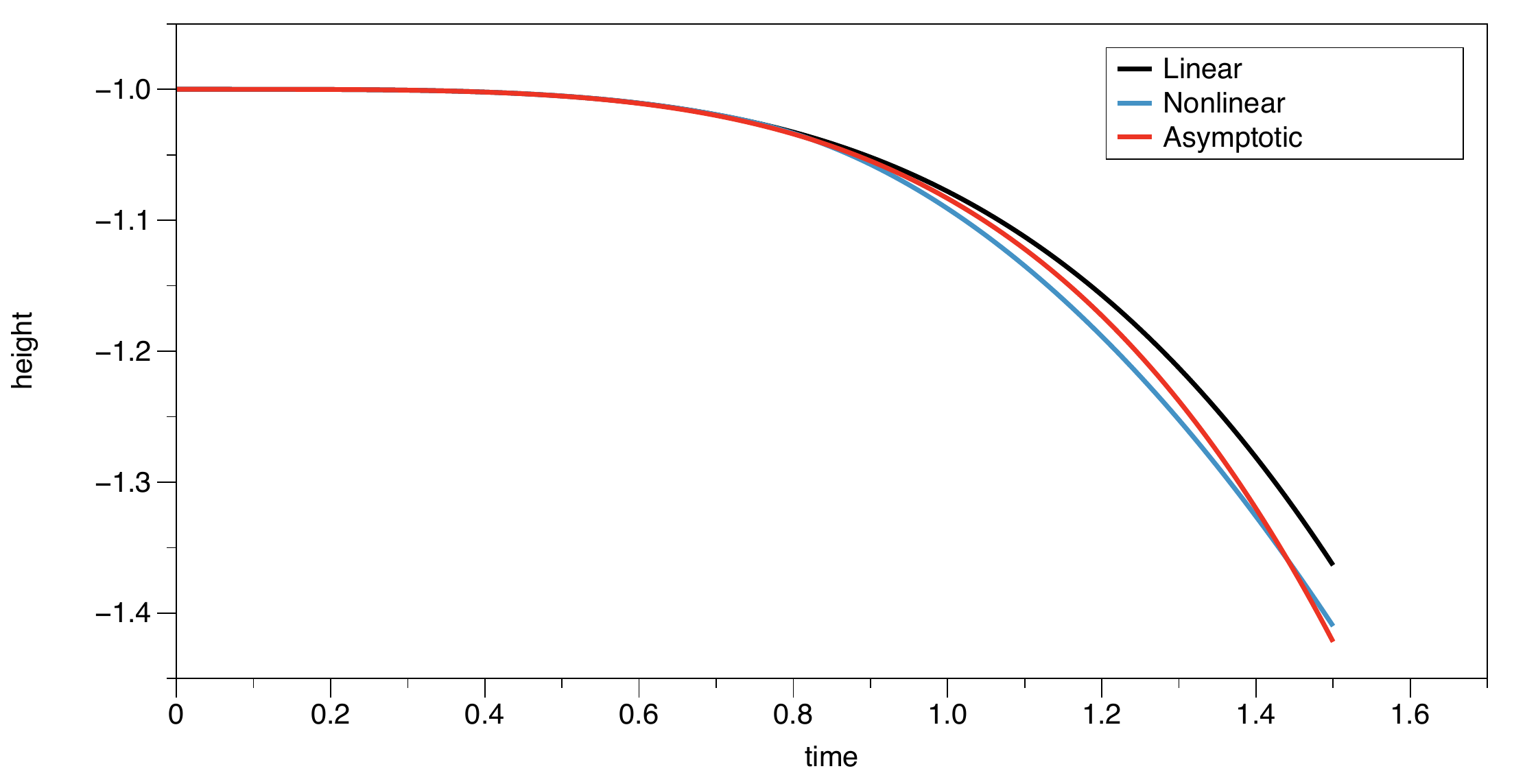}
\caption{Left panel:  experimental, asymptotics, and full linear solution comparison for the trial depicted in Figure 1 for mass, $m_4$ in cgs units.  Right panel:  Comparison of asymptotic, full linear solution, and full nonlinear solution for case with two equal masses and one spring, in non-dimensional coordinates.   }
\label{full}
\end{figure}

Lastly, it is interesting to compare the complete dynamics of the model linear mass-spring system with the experimental trajectories and the short time asymptotics.  It is clear from observing the montage in Figure \ref{montage} that once the top spring collapses upon mass, $m1$, the collision is not elastic, and the spring essentially stops providing any internal dynamics for the system except as to continue to provide inertia from its mass.  A perhaps better model would be one which includes nonlinear cut-offs to remove a spring once its compressed length is reduced beyond some threshold.  To this end, we study more in depth the bottom mass track of the first trials with its short time asymptotics as well as with the full solution of the coupled linear system.  This comparison is presented in the left panel of Figure \ref{full}.  As observed before, the short time asymptotic prediction is excellent here.  One interesting thing to note is that the full solution to the linear system appears to fall at a slower rate than both the data and the asymptotic theory beyond the first tenth of a second, with the data perhaps on those timescales predicted better by the asymptotic theory than the full linear solution.  This is likely on account of the springs role in driving internal dynamics ceasing after the collisions.

To demonstrate and study this, we can incorporate a nonlinearity into the system in the form of a cut-off which turns the spring force off once the spring has recoiled to its un-stretched position.  To this end, it is useful to solve this system in a piecewise solution for two identical masses, connected by one spring using the integral of the motion given by the center of mass:  In dimensional coordinates:
\begin{eqnarray*}
m\frac{d^2 x}{dt^2} &=& - k (x - y) f(x-y) -m g\\
m\frac{d^2 x}{dt^2} &=&  k (x - y) f(x-y) -m g\\
x(0)&=&0, \quad y(0)=\frac{-mg}{k}, \quad x'(0)=y'(0)=0
\end{eqnarray*}
where the cut-off, $f(x)$ is a Heaviside function taking zero value if $x<L_s$, where $L_s$ is the spring length.  This system can be decoupled through adding the two odes, yielding $y(t)=-t^2-1-x(t)$.  For short time, the nonlinearity is identically equal to one, and that ode for $x$ is solved directly as $x(t)=-\frac{mg}{2k} -\frac{g t^2}{2} + \frac{mg}{2k} \cos{\sqrt{\frac{2k}{m}}t}$.  This solution is valid until the spring relaxes with $x(t^*)=L_s$ at which point both $x$ and $y$ switch to pure free fall, initialized by the positions and speeds at $t=t^*$.  The comparison of these two evolutions for $y(t)$ are presented in the right panel of Figure \ref{full} in non-dimensional coordinates, here taking $L_s=1/2$. Similar to the experimental case reported in the left panel, the full nonlinear solution falls more quickly on the depicted timescales than the full linear solution.  Also shown is the short time asymptotic solution, $-1-t^4/12$, for this system in non-dimensional coordinates.  We remark that similar analyses can be performed for the complete system, but as this does not impact our original goal of developing a theory for the hang-time, we will pursue this in future work.

\section{Discussion}

Here, we have studied experiments with dropping a hanging mass-spring system, initially in equilibrium under tension.  These experiments show a delayed hang-time for the bottom mass:  the top mass is released and falls before any of the other masses below start moving.   The springs are measured to be linear Hooke law springs, and a complete theory for each masses' short time asymptotics is developed using a combination of Cramer's rule and Tauberian identities, which shows the $j'th$ mass height falls as $C_j t^{2j}$.  This theory is confirmed on two trials, one with no explicit top mass, and a second with a top mass attached.  This latter case shows a substantially enhanced hang-time as compared with the former case.  Further, long time asymptotics are derived which show the long time solutions fall with the acceleration of gravity, with corrections being subdominant.  Lastly, the full linear dynamics are compared with the experiment, which shows the trend that the linear dynamics fall slower than the experiment.  A nonlinear model with a force cut-off is presented and solved for the two body case, which depicts the same trend.  The success of the linear theory is somewhat expected as everything is started from rest in equilibrium.  This points to the expectation that a general non-linear system with non-linear spring laws presumably experiences the same short time asymptotics, assuming the spring law possesses a non-zero linear term.  

A few comments are in order.  As discussed above, the short time asymptotics are power law, which mathematically has the property that all the masses are moving for any positive time.  The theory predicts that the $j$'th mass falls as 
$Q_j \frac{t^{2j} }{(2j)!}$.  While each of these functions is non-zero at short time, the motion will not be measurable until a threshold (such as one pixel in a camera image) is passed.  Moreover, the solutions further down the chain are flatter and flatter.  This is somewhat reminiscent of the heat equation which has the peculiar property that a compactly supported initial data gives rise to a non-zero solution for any positive time, {\it arbitrarily far from} the region of the compact data.  Despite violating all sorts of physical principles, this behavior is immediately rectified by observing that no information is propagated until the solution reaches a measurable threshold.  As such, the heat equation does not propagate information at infinite speed, but rather following the usual diffusive scaling laws in which information propagates at sub-ballistic (diffusive) rates.  We remark that this modeling peculiarity of the heat equation could be further rectified directly by use of a Telegrapher's equation for the heat transport, which is wavelike at short times.  An interesting derivation of Telegrapher's equation in the context of random walks may be found in work of Goldstein \cite{Goldstein,Taylor} .  A similar approach can be followed here in which a continuum limit of our mass-spring system gives rise to a hyperbolic wave equation.  Such equations enjoy ballistic time propagation of information.  It is an interesting question to examine the continuum limit of our system, particularly as a first model for a dropped slinky which shows a similar delay timescale.  We leave that to future work.

Of course, Wile E. Coyote will always be tricked into falling off cliffs by the clever roadrunner.  It is our hope that perhaps Wile E. may learn from the present discussion and hopefully apply our findings to situations involving long, hanging ropes, cut but under tension which our calculations provide some indications of the timescales under which the Coyote may have time to climb said rope without falling to his eventual, but temporary, demise.  

We feel that this example, in addition to possessing interesting complex internal dynamics, provides exciting new curriculum content for modern advanced mathematical training in the sciences as it shows unexpected behavior and further provides a data rich experimental system which can be used to demonstrate the power of mathematical modeling using somewhat unfamiliar asymptotic methods and complex analysis.  

\section{Acknowledgements}

We would like to thank Kevin Guskiewicz who was Dean of the UNC College of Arts and Sciences, and Chris Clemens who was Senior Associate Dean for the Sciences for supporting our interdisciplinary honors differential equations class for several years in which this system and others were explored.  We thank David Adalsteinsson for his wonderful Image Tank and Data Graph programs.  We thank David McLaughlin for helpful discussions and thank Jim Mahaney for help with the experiments.  We thank funding through the National Science Foundation Grant Nos.:DMS-1910824; and the Office of Naval Research Grant No: ONR N00014-18-1-2490. 

\bibliographystyle{abbrv}

\end{document}